\documentclass{gtart_h}  

\def\ifplaintex{\expandafter\ifx\csname documentclass\endcsname\relax}

\ifplaintex 
\else
\fi

\expandafter\ifx\csname beginpicture\endcsname\relax
\expandafter\ifx\csname documentclass\endcsname\relax
\input pictex \else
\input prepictex \input pictex \input postpictex \fi\fi

\def\gt{{\mathsurround=0pt\it $\cal G\mskip-2mu$eometry \&\ 
$\cal T\!\!$opology}}        %  journal title in recommended style

\def\gtp{{\mathsurround=0pt\it $\cal G\mskip-2mu$eometry \&\ 
$\cal T\!\!$opology $\cal P\!$ublications}}  % GT publications

\def\lognumber#1{\def\thelognumber{#1}}
\def\volumenumber#1{\def\thevolumenumber{#1}}
\def\papernumber#1{\def\thepapernumber{#1}}
\def\volumeyear#1{\def\thevolumeyear{#1}}

\def\pagenumbers#1#2{\def\startpage{#1}\def\finishpage{#2}}
\def\published#1{\def\publishdate{#1}}
\def\proposed#1{\def\theproposer{#1}}
\def\seconded#1{\def\theseconders{#1}}
\def\received#1{\def\receiveddate{#1}}
\def\revised#1{\def\reviseddate{#1}}
\def\accepted#1{\def\accepteddate{#1}}

\def\asciiaddress#1{\def\theasciiaddress{#1}}
\def\asciiemail#1{\def\theasciiemail{#1}}

\long\def\asciiabstract#1{\long\def\theasciiabstract{#1}}

\def\shorttitle#1{\def\theshorttitle{#1}}

\let\\\par\let\thelognumber\relax
\let\thevolumenumber\relax\let\thepapernumber\relax
\let\thevolumeyear\relax\let\thesamplenumber\relax\let\startpage\relax
\let\finishpage\relax\let\publishdate\relax\let\receiveddate\relax
\let\reviseddate\relax\let\accepteddate\relax\let\theasciititle\relax
\let\theasciiauthors\relax\let\theasciiaddress\relax
\let\theasciiabstract\relax
\let\theasciiemail\relax\let\theshortauthors\relax\let\theshorttitle\relax

\long\def\maketitlep{   % start of definition of \maketitlep

\count0=\startpage

\gt\hfill      %   Journal title (top left) 
\beginpicture
\setcoordinatesystem units <0.33truein, 0.33truein> point at 2.2 0.9
\setplotsymbol ({$\cal G$})
\plotsymbolspacing=9truept
\circulararc 315 degrees from 0 1 center at 0 0
\setplotsymbol ({$\cal T$})
\circulararc 315 degrees from 1 -1 center at 1 0
\endpicture
\break
{\small\ifx\thesamplenumber\relax % sample?  
Volume \else Sample
\fi\thevolumenumber\ (\thevolumeyear)
\startpage--\finishpage\nl
Published: \publishdate}
\vglue 0.5truein plus 0.4fil minus 0.1truein

{\parskip=0pt\leftskip 0pt plus 1fil\def\\{\par\smallskip}{\ifplaintex\large
\else\Large\fi\bf\thetitle}\par\medskip}   

\vglue 0pt plus 0.1fil 

{\parskip=0pt\leftskip 0pt plus 1fil\def\\{\par}{\sc\theauthors}
\par\medskip}

\vglue 0pt plus 0.1fil 

{\small\parskip=0pt\let\newline\\
{\leftskip 0pt plus 1fil\def\\{\par}{\sl\theaddress}\par}
\expandafter\ifx\theemail\relax    % email address?
\relax\else\vglue 5pt plus 0.02fil minus 2pt\def\\{\stdspace{\rm 
and}\stdspace} 
\cl{Email:\stdspace\tt\theemail}\fi
\ifx\theurl\relax                  % URL given?
\relax\else\vglue 5pt plus 0.02fil minus 2pt\def\\{\stdspace{\rm 
and}\stdspace}
\cl{URL:\stdspace\tt\theurl}\fi\par}

\vglue 7pt plus 0.3fil minus 3pt

{\bf Abstract}
\vglue 5pt plus 0.1fil minus 2pt

\theabstract

\vglue 7pt plus 0.3fil minus 3pt

{\bf AMS Classification numbers}\quad Primary:\quad \theprimaryclass

Secondary:\quad \thesecondaryclass

\vglue 5pt plus 0.3fil minus 2pt

{\bf Keywords:}\quad \thekeywords

\vglue 10pt plus 0.5fil minus 5pt

{\small  Proposed: \theproposer\hfill Received: \receiveddate\nl
Seconded: \theseconders\hfill 
\ifx\reviseddate\relax                         % paper revised?
Accepted: \accepteddate                        % no
\else
Revised: \reviseddate                          % yes
\fi}
\eject
}       %  end of definition of \maketitlep

\let\maketitlepage\maketitlep
\let\maketitle\maketitlepage

\font\phead=cmsl9 scaled 950
\font=cmsl9 scaled 1050
\font\pnum=cmbx10 scaled 913
\font\lnum=cmbx10 
\font\pfoot=cmsl9 scaled 950
\font=cmsl9 scaled 1050
\ifplaintex
\headline{\vbox to 0pt{\vskip -4.5mm\line{\small\phead\ifnum
\count0=\startpage ISSN 1364-0380 (on line)
1465-3060 (printed) \hfill {\pnum\folio}\else\ifodd\count0\def\\{ }% 
\ifx\theshorttitle\relax\thetitle\else\theshorttitle\fi\hfill{\pnum\folio}
\else\def\\{ and }{\pnum\folio}\hfill\ifx\theshortauthors\relax\theauthors
\else\theshortauthors\fi\fi\fi}\vss}}
\footline{\vbox to 0pt{\vglue 0mm\line{\small\pfoot\ifnum\count0=\startpage
\copyright\ \gtp\hfill\else
\gt, Volume \thevolumenumber\ (\thevolumeyear)\hfill\fi}\vss
}}
\else
\makeatletter
\def\@oddhead{{\small\lhead\ifnum\count0=\startpage ISSN 1364-0380 (on line)
1465-3060 (printed) \hfill {\lnum\number\count0}\else\ifodd\count0
\def\\{ }\ifx\theshorttitle\relax \thetitle \else\theshorttitle\fi\hfill
{\lnum\number\count0}\else\def\\{ and }{\lnum\number\count0}
\hfill\ifx\theshortauthors\relax 
\theauthors\else\theshortauthors\fi\fi\fi}}\def\@evenhead{\@oddhead}
\def\@oddfoot{\small\lfoot\ifnum\count0=\startpage\copyright\ \gtp\hfill\else
\gt, Volume \thevolumenumber\ (\thevolumeyear)\hfill\fi}
\def\@evenfoot{\@oddfoot}
\makeatother
\fi

\newwrite\gtoutfile
\long\gdef\makeheadfile{  %%% start of definition of \makeheadfile
{\def\\{, }\def\s{ }
\immediate\openout\gtoutfile head.xxx
\immediate\write\gtoutfile{Proxy-for: \ifx\theasciiauthors\relax
\theauthors\else\theasciiauthors\fi\s<\ifx\theasciiemail\relax\theemail\else\theasciiemail\fi>}
\immediate\write\gtoutfile{\noexpand\\}
\immediate\write\gtoutfile{Authors: \ifx\theasciiauthors\relax
\theauthors\else\theasciiauthors\fi}
{\def\\{ }\immediate\write\gtoutfile{Title: \ifx\theasciititle\relax
\thetitle\else\theasciititle\fi}}
\immediate\write\gtoutfile{Subj-class: GT or SG or MG etc}
\immediate\write\gtoutfile{MSC-class: \theprimaryclass\ifx\thesecondaryclass\relax\else, \thesecondaryclass\fi}
\immediate\write\gtoutfile{Journal-ref: Geom. Topol. \thevolumenumber
(\thevolumeyear) \startpage-\finishpage}
\immediate\write\gtoutfile{Comments: Published by Geometry and Topology at}
\immediate\write\gtoutfile{\s\s http://www.maths.warwick.ac.uk/gt/GTVol\thevolumenumber/paper\thepapernumber.abs.html}
\immediate\write\gtoutfile{\noexpand\\}
\immediate\write\gtoutfile{}
\ifx\theasciiabstract\relax
\immediate\write\gtoutfile{\theabstract}\else
\immediate\write\gtoutfile{\theasciiabstract}\fi
\immediate\write\gtoutfile{}
\immediate\write\gtoutfile{\noexpand\\}
\immediate\write\gtoutfile{}
\immediate\closeout\gtoutfile}}  %%% end of definition of \makeheadfile

\def\maketitlepage{\maketitlep\makeheadfile}
\let\maketitle\maketitlepage

\lognumber{368}
\volumenumber{8}\papernumber{13}\volumeyear{2004}
\pagenumbers{539}{564}
\received{9 October 2003}
\revised{16 March 2004}
\published{17 March 2004}
\accepted{23 December 2004}
\proposed{Jean-Pierre Otal}
\seconded{David Gabai, Martin Bridson}

\usepackage{amsmath, amssymb,
graphics}

\def\S{Section }
\def\leq{\leqslant}
\def\geq{\geqslant}
\newtheorem{theorem}{Theorem}
\newtheorem{proposition}[theorem]
{Proposition}
\newtheorem{lemma}[theorem]{Lemma}

\newtheorem{corollary}[theorem]{Corollary}

\theoremstyle{remark}
\newtheorem{remark}[theorem]{Remark}

\begin{document}

\title {The metric
space of geodesic\\laminations on a surface: I}
\shorttitle{The metric space of
geodesic laminations}

\author{Xiaodong Zhu\\Francis Bonahon}

\address {NetScreen
Technologies,  805 11-th Avenue\\Building 3, 
Sunnyvale, CA 94089,
USA\\{\rm and}\\Department
of Mathematics,  University of
Southern California\\Los Angeles,
CA 90089-1113, USA}

\asciiaddress{NetScreen
Technologies,  805 11-th Avenue\\Building 3, 
Sunnyvale, CA 94089,
USA\\and\\Department
of Mathematics,  University of
Southern California,\\Los Angeles,
CA 90089-1113, USA}

\gtemail{\mailto{xzhu@netscreen.com}{\rm\qua and\qua}\mailto{fbonahon@math.usc.edu}}
\asciiemail{xzhu@netscreen.com, fbonahon@math.usc.edu}

\urladdr{http://math.usc.edu/~fbonahon}

\begin{abstract} We consider the
space of geodesic laminations on
a surface, endowed with the
Hausdorff metric $d_{\mathrm H}$
and with a variation of this
metric called the $d_{\log}$
metric. We compute and/or
estimate the Hausdorff dimensions
of these two metrics. We also
relate these two metrics to
another metric which is
combinatorially defined in terms
of train tracks.
\end{abstract}

\asciiabstract{We consider the space of geodesic laminations on a
surface, endowed with the Hausdorff metric d_H and with a variation of
this metric called the d_log metric. We compute and/or estimate the
Hausdorff dimensions of these two metrics. We also relate these two
metrics to another metric which is combinatorially defined in terms of
train tracks.}


\maketitlepage

Let $S$ be a compact surface,
possibly with boundary. We
restrict attention to the case
where the Euler characteristic of
$S$ is negative.

Geodesic laminations on $S$ were
introduced by Bill Thurston
\cite{Thu76, Thu84} to
provide a completion for the
space of simple closed curves on
$S$. They also occur in various
problems in low-dimensional
topology and geometry, and have
been used as a very successful
tool to attack these problems. So
far, most of their use and
analysis has been focused on
\emph{measured} geodesic
laminations, namely geodesic
laminations endowed with the
additional structure of a
transverse measure.  However,
there are several contexts where
geodesic laminations occur
without a preferred transverse
measure. A fundamental example of
this is the ending lamination of a
geometrically infinite hyperbolic
3--manifold; see for instance
\cite{Min02} for a discussion, and
compare
\cite{MasMin99, MasMin00, Kla}.
The current article and its
sequel \cite{BonZhu3} find
their motivation in an attempt to
systematically  develop an
`unmeasured' theory of geodesic
laminations. In particular, one
could hope that the space
$\mathcal L(S)$ of all geodesic
laminations on $S$ would share
many features with the space
$\mathcal {PML}(S)$ of projective
measured geodesic laminations,
which was analyzed in
\cite{Thu76,
Thu84, FatLauPoe, PenHar} and
in subsequent work. To some
extent, our results show that
this is not the case.

We endow  the
space
$\mathcal L(S)$ of geodesic
laminations on $S$ with the
Hausdorff metric
$d_{\mathrm H}$, whose definition
can be found in
Section~\ref{sect:GeodLam}. Our
first result is the following.

\begin{theorem}
\label{thm:DimZero} The
Hausdorff dimension of the metric
space
$\left( \mathcal L(S), d_{\mathrm
H}\right)$ is equal to~$0$.
\end{theorem}

	This should not be confused
with the result of Joan Birman and
Caroline Series
\cite{BirSer} which says that the
union of all  geodesic
laminations forms a subset of $S$
of  Hausdorff dimension 1,
although the method of  proof is
closely related to the main
ingredient of  \cite{BirSer}.
Thurston already proved in
\cite[\S 10]{Thu86} that the
subset
$\mathcal L^{\mathrm{cr}} (S)$ of
$\mathcal L(S)$ that is the
closure of the set of
multicurves (= disjoint
unions of finitely many
simple closed
geodesics) has Hausdorff dimension
0. Note that
$\mathcal L(S)$ is significantly
larger than
$\mathcal L^{\mathrm{cr}} (S)$,
as can for instance be seen in the
examples of
\cite{BonZhu3}.

An immediate corollary of
Theorem~\ref{thm:DimZero} is that
the space $\mathcal L(S)$ also
has topological dimension 0 (see
for instance
\cite[Proposition 2.5]{Fal})

\begin{corollary} The space
$\mathcal L(S)$ is totally
disconnected.
\end{corollary}

Our analysis \cite{BonZhu3} of the
simple examples where $S$ is the
once punctured torus or the
four-times punctured sphere shows
the kind of intricacies which one
can expect for the topology of
$\mathcal L(S)$. In these
examples, $\mathcal L(S)$
consists of a standard Cantor set
$K$ in $\mathbb R$ with finitely
many isolated points in each of
the gaps of $\mathbb R-K$.

Theorem~\ref{thm:DimZero} is at
first somewhat disappointing, as
one would prefer a dimension which
reflects the topological
complexity of $S$. On second
thoughts, its proof suggests a
different metric on
$\mathcal L(S)$ which has better
invariance properties than the
Hausdorff metric
$d_{\mathrm H}$ but induces the
same topology. Indeed, the
definition of $d_{\mathrm H}$
requires that we choose a
negatively curved metric on $S$,
and a different choice for this
metric changes $d_{\mathrm H}$ by
a H\"older equivalence. In
particular, the only way the
Hausdorff dimension of
$d_{\mathrm H}$ could be
independent of this choice is if
the dimension is 0 or $\infty$,
which in retrospect makes
Theorem~\ref{thm:DimZero} much
more predictable.

For this reason, it is much
better to consider the metric
$d_{\log}$ introduced in
Section~\ref{sect:dlog}. For small
distances, $d_{\log}$ is just
$1/\left\vert \log d_{\mathrm
H}\right\vert$. The metric
$d_{\log}$ depends only on the
topology of $S$ up to Lipschitz
equivalence. In particular, its
Hausdorff
dimension depends only on the
topology of $S$. Note that
$d_{\log}$ and $d_{\mathrm H}$
induce the same topology on
$\mathcal L(S)$.

A similar idea of
modifying the metric, or at
least of using a different gauge
to compute Hausdorff dimensions,
was already used by Geoff Mess in
\cite[\S5]{Mes}, in the context of
\cite{BirSer}.

\begin{theorem}
\label{thm:HausDimL(S)1}
If the
surface
$S$ is different from the three
times punctured sphere, the twice
punctured projective plane or the
once punctured Klein bottle, the
Hausdorff dimension $\left(
\mathcal L(S),  d_{\log} \right)$
is finite, and at least equal to
$2$.
\end{theorem}

In the three exceptional cases
where $S$ is the three times
punctured sphere, the twice
punctured projective plane or the
once punctured Klein bottle, the
space
$\mathcal L(S)$ is finite, and
consequently has Hausdorff
dimension 0; see for instance
\cite[Proposition 16]{BonZhu3}.

The lower estimate comes from
\cite{BonZhu3}, where we prove
that the Hausdorff dimension of
$\left(
\mathcal L(S),  d_{\log} \right)$
is exactly equal to 2 when $S$ is
the once punctured torus or the
4-times punctured sphere.
Proposition~\ref
{thm:dimRoughBounds}
provides an explicit upper bound
for the Hausdorff dimension, but
this upper bound is far from being
sharp. By analogy with the
arguments of \cite{BonZhu3}, we
conjecture that the Hausdorff
dimension
$\delta$ of
$\left(
\mathcal L(S),  d_{\log} \right)$
is always an integer (although
the value of the conjectured
integer tends to depend on which
heuristic argument we are using).
However, in the cases of the once
punctured  torus and the 4-times
punctured sphere the
$\delta$--Hausdorff measure of
$\left(
\mathcal L(S),  d_{\log} \right)$
is the zero measure, and we
conjecture that this is always
the case. To some extent, this is
another disappointment as it does
not yield as interesting a
dynamical systems as one could
have hoped for.

In the last part of the paper,
we show that the metric
$d_{\log}$ has a strong
combinatorial flavor, in terms
of train tracks on the surface
$S$. More precisely, the geodesic
lamination
$\lambda$ is \emph{weakly carried}
by the train track $\Theta$ if,
for every leaf $l$ of $\lambda$,
there is a bi-infinite curve $c$
which is immersed in
$\Theta$ and is homotopic to $l$ by
a homotopy which moves points by
a bounded amount. In this case, we
say that the leaf $l$ is
\emph{tracked} by the curve $c$.
When $c$ crosses a succession of
edges $e_1$, $e_2$, \dots, $e_n$
of $\Theta$, in this order, we say
that $\lambda$ \emph{realizes} the
edge path $\left\langle e_1, e_2,
\dots, e_n \right\rangle$.

Given two geodesic laminations
$\lambda$, $\lambda'$ carried by
the train track $\Theta$, we define
their \emph{combinatorial
distance} as follows.
\begin{equation*}
d_\Theta(\lambda,
\lambda'){=}\min\! \left\{ \frac1{r+1};
\lambda \text{ and } \lambda'
\text{ realize the same edge
paths of length } r\text{ in }\Theta
\right\}
\end{equation*}

The sets $\mathcal L^{\mathrm
w}(\Theta)$, consisting of those
geodesic laminations which are
weakly carried by the train track
$\Theta$, can be used as local
charts for the space $\mathcal
L(S)$ of all geodesic
laminations, in the sense that
$\mathcal L(S)$ can be written as
the union of finitely many $\mathcal L^{\mathrm
w}(\Theta_i)$.

We prove the following estimate,
which relates the Hausdorff
metric $d_{\mathrm H}$ to the
combinatorial metric $d_\Theta$.

\begin{theorem}
\label{thm:HausCombEstimate1}
  Let $\Theta$ be a
train track on the
negatively curved surface
$S$. Then, there exists two
constants $c_1$, $c_2>0$
such that, for every geodesic
laminations
$\lambda$, $\lambda'$ which are
weakly carried by $\Theta$,
\begin{equation*}
\mathrm e^{-c_1/d_{\Theta}
\left( \lambda,
\lambda' \right)}
\leq
d_{\mathrm H} \left( \lambda,
\lambda' \right)
\leq  \mathrm e^{-c_2 /d_{\Theta}
\left( \lambda,
\lambda' \right)}
\end{equation*}
\end{theorem}

The first inequality is harder to
prove, and deeper, than
the second one. Indeed, this
first inequality asserts in
particular that two geodesic
laminations which are close for
the Hausdorff topology must be
combinatorially close.

Since, for small distances, the
metric $d_{\log}$ is just
$1/\left|\log d_{\mathrm
H}\right|$, the above estimate
can be rephrased in the
following way.

\begin{corollary}
There exists two
constants $c$, $c'>0$
such that, for every pair of geodesic
laminations
$\lambda$, $\lambda'$ which are
weakly carried by $\Theta$,
\begin{equation*}
cd_{\Theta}
\left( \lambda,
\lambda' \right)
\leq
d_{\log} \left( \lambda,
\lambda' \right)
\leq  c'd_{\Theta}
\left( \lambda,
\lambda' \right).
\end{equation*}
In other words,  the restrictions
of the metrics $d_{\log}$ and
$d_\Theta$ to the subspace of
$\mathcal L(S)$ consisting of
those geodesic laminations which
are weakly carried by the train
track are
Lipschitz equivalent.
\end{corollary}

Theorem~\ref
{thm:HausCombEstimate1} is a key
ingredient for the arguments of
\cite{BonZhu3}, which themselves
provide the lower estimate of
Theorem~\ref
{thm:HausDimL(S)1}.

The results of this article were
first proved in the dissertation
\cite{Zhu}. The authors are
grateful to the referee for
several valuable comments on
the first version of this
article.

{This work was partially
supported by grants DMS-9504282,
DMS-9803445 and DMS-0103511 from
the National Science Foundation.}
\date{\today}

\section{Geodesic laminations}
\label{sect:GeodLam}

Let $S$ be a compact surface,
possibly with boundary, whose
Euler characteristic
$\chi(S)$ is negative. The
hypothesis on the Euler
characteristic implies that
$S$ admits a riemannian metric
$m$ of negative curvature%
\footnote{The metric can even be
chosen to have constant
curvature $-1$. Those readers who
are more comfortable with
hyperbolic geometry can therefore
assume that this is the case,
although this does not really
simplify the arguments.} for
which the boundary
$\partial S$ is totally geodesic.
Fix such a metric.

A \emph{geodesic lamination} is a
non-empty closed subset of $S$
which can be decomposed as the
union of a family of disjoint
simple complete geodesics. Recall
that a complete geodesic is one
which admits a parametrization by
arc length defined over all of
$\left]-\infty, +\infty\right[$.
In particular, a complete
geodesic does not transversely
hit the boundary. A closed
geodesic is always complete, and
a geodesic is allowed to be
contained in the boundary
$\partial S$. A geodesic is
\emph{simple} it it has no
transverse self-intersection. It
can be shown that the
decomposition of a geodesic
lamination into disjoint simple
geodesics is unique. These
geodesics are called the
\emph{leaves} of the geodesic
lamination. We refer to
\cite{CanEpsGre, CasBle,
PenHar, Bon10} for general
background on geodesic
laminations.

A fundamental example of geodesic
lamination is provided by the
union of finitely many disjoint
simple closed geodesics. A
geodesic lamination of this type
is called a \emph{multicurve}.
However, typical geodesic
laminations are much more complex
from a topological and dynamical
point of view.

This paper is devoted to studying
the space $\mathcal L(S)$ of all
geodesic laminations of $S$.

To define a metric on $\mathcal
L(S)$, we will use the classical
\emph{Hausdorff metric} on the set
$C(X)$ of all non-empty bounded
closed subsets of a metric space
$(X,d)$, defined by
\begin{equation*} d_{\mathrm
H}(C,C') =
\inf\left\{\varepsilon; C\subset
N_\varepsilon(C')
\text{ and } C'\subset
N_\varepsilon(C)
\right\}\quad \text{if}\quad
C, C'\in
C(X),
\end{equation*} 
where
$N_\varepsilon(C)$ denotes the
$\varepsilon$--neighborhood
$\left\{x\in X;
\exists y\in C,
d(x,y)\leq\varepsilon
\right\}$  of $C$. It is
well-known that
$d_{\mathrm H}$ defines a metric
on
$C(X)$; see for instance
\cite{Mun}.

To define a metric on $\mathcal
L(S)$, we could consider it as a
subset of the set $C(S)$ of all
closed non-empty subsets of $S$,
and use the restriction of the
Hausdorff metric $d_{\mathrm H}$.
However, for reasons explained
later in this section, it makes
more sense to use the
\emph{projective tangent bundle}
$PT(S)$ of $S$, defined as the
set of pairs $(x,l)$ where $x\in
S$ and where
$l$ is a line (through the
origin) in the tangent space
$T_xS$. A geodesic lamination
lifts to a closed subset
$\widehat\lambda$ of $PT(S)$ by
considering all $(x,T_x\lambda)$
where
$x\in\lambda$ and
$T_x\lambda\subset T_xS$ is the
line tangent to the leaf of
$\lambda$ passing through $x$.

In this way, we have identified a
geodesic lamination
$\lambda\in\mathcal L(S)$ to a
closed subset $\widehat\lambda$
of $PT(S)$. The Levi-Civita
connection of the riemannian
metric $m$ on
$S$ enables us to lift $m$ to a
riemannian metric on $PT(S)$,
thereby turning $PT(S)$ into a
metric space. We then define the
\emph{Hausdorff distance} between
the geodesic laminations
$\lambda$ and
$\lambda'\in\mathcal L(S)$ as the
Hausdorff distance $d_{\mathrm
H}(\widehat\lambda,
\widehat\lambda')$ between the
closed subsets
$\widehat\lambda$,
$\widehat\lambda'\in C(PT(S))$.
We will write
$d_{\mathrm H}(\lambda, \lambda')$
for $d_{\mathrm
H}(\widehat\lambda,
\widehat\lambda')$.

This definition of $d_{\mathrm
H}(\lambda, \lambda')$ can also
be rephrased as
\begin{equation*}
d_{\mathrm H}(\lambda, \lambda')=
\min
\left\{\varepsilon; \,
\begin{aligned}
\forall x\in\lambda, \exists
x'\in \lambda' ,d\left((x,
T_x\lambda), (x',
T_{x'}\lambda')\right)
<\varepsilon\\
\forall x'\in\lambda', \exists
x\in \lambda, d\left((x,
T_x\lambda), (x',
T_{x'}\lambda')\right)
<\varepsilon
\end{aligned}
\right\}
\end{equation*}
where $d$ here denotes the
metric in the projective tangent
bundle
$PT(S)$.

\begin{lemma}[{\cite[\S
3]{CasBle}},
{\cite[\S4.1]{CanEpsGre}},
{\cite[\S 1.2]{Bon10}}]
\label{thm:L(S)compact} The
metric space
$\left( \mathcal L(S), d_{\mathrm
H} \right)$ is compact.\qed
\end{lemma}

We now investigate what happens
when we replace the negatively
curved metric $m$ on $S$ by
another negatively curved metric
$m'$ with totally geodesic
boundary. By compactness of $S$,
the distortion between $m$ and
$m'$ is bounded, and a complete
$m$--geodesic
$g$ is quasi-geodesic for $m'$; a
general property of
quasi-geodesics in negatively
curved manifolds then says that
there exists a unique
$m'$--geodesic $g'$ which is
homotopic to $g$ by a homotopy
moving points by a bounded
amount.  This correspondence
between $m$--geodesics and
$m'$--geodesics sends simple
geodesics to simple geodesics,
and disjoint geodesics to
disjoint geodesics, and therefore
induces a one-to-one
correspondence between the set
$\mathcal L(S,m)$ of $m$--geodesic
laminations and the set $\mathcal
L(S,m')$ of $m'$--geodesic
laminations. In other words, the
set $\mathcal L(S,m)$ is really
independent of the metric $m$,
which justifies our original
notation $\mathcal L(S)=\mathcal
L(S,m)$.

Let us see how the Hausdorff
metric
$d_{\mathrm H}$ depends on the
metric
$m$.  Two metrics $d$ and $d'$ on
a set
$X$ are said to be \emph{H\"older
equivalent} if there exists
constants $C_1$, $C_2$, $\nu_1$,
$\nu_2>0$ such that
\begin{equation*}
C_1 d(x,y)^{\nu_1} \leq
d'(x,y)\leq C_2 d(x,y)^{\nu_2}
\end{equation*}
for every $x$,
$y\in X$. They are
\emph{Lipschitz equivalent} if
they are H\"older equivalent with
exponents
$\nu_1=\nu_2=1$.

\begin{lemma}
\label{lem:HolderInv} As we
change the negatively curved
metric $m$ on $S$, the Hausdorff
metric
$d_{\mathrm H}$ of
$\mathcal L(S)$ varies only by a
H\"older equivalence.
\end{lemma}

\begin{proof} Let $m$ and $m'$ be
two negatively curved metrics on
$S$ for which the boundary is
totally geodesic.

In general, there is no
homeomorphism
$S\rightarrow S$ which sends each
$m$--geodesic $g$ to the
corresponding
$m'$--geodesic. However, there is
always a homeomorphism
$\varphi\co PT(S)\rightarrow
PT(S)$ which sends the lift
$\widehat g\subset PT(S)$ of each
$m$--geodesic $g$ to the lift
$\widehat g'\subset PT(S)$ of the
corresponding
$m'$--geodesic $g'$. See
\cite[\S 7.2]{Gro} or \cite[\S
19.1]{KatHas}. In addition,
$\varphi$ is \emph{H\"older
bicontinuous} in the sense that
there exists constants $C_1$, $C_2$, $\nu_1$,
$\nu_2>0$ such that
\begin{equation*}
C_1 d(x,y)^{\nu_1} \leq
d'(\varphi(x), \varphi(y))\leq C_2
d(x,y)^{\nu_2}
\end{equation*} for every $x$,
$y\in PT(S)$, where $d$ and $'$
are the metrics on $PT(S)$
respectively defined by $m$ and
$m'$. It immediately follows that
the Hausdorff metrics $d_{\mathrm
H}$ and
$d_{\mathrm H}'$ defined by $m$
and
$m'$ on
$\mathcal L(S)  =\mathcal L(S,m)
=\mathcal L(S,m')$ are H\"older
equivalent.
\end{proof}

As indicated earlier, instead of
$d_{\mathrm H}$, we could have
considered on $\mathcal L(S)$ the
restriction $\delta_{\mathrm H}$
of the Hausdorff metric of the
spaces of closed subsets of $S$
(as opposed to
$PT(S)$). The identity map
$\left(\mathcal L(S), d_{\mathrm
H}\right)
\rightarrow
\left(\mathcal L(S),
\delta_{\mathrm H}\right)$ is
continuous by definition of the
topologies, and is therefore a
homeomorphism by the compactness
property of
Lemma~\ref{thm:L(S)compact}. In
other words,
$\delta_{\mathrm H}$ induces the
same topology as
$d_{\mathrm H}$ on $\mathcal
L(S)$ and, consequently, is
exactly as good as $d_{\mathrm
H}$ from a topological point of
view.  However, from a metric
point of view,
Lemma~\ref{lem:HolderInv} does
not seem to hold for
$\delta_{\mathrm H}$ in any easy
way. In addition to this technical
problem, it may philosophically
make more sense to consider a
geodesic lamination as a subset
of $PT(S)$, as indicated by the
general framework of geodesic
currents developed in
\cite{Bon86, Bon88}.

\section{Fattened train tracks}
\label{sect:FatTrainTracks}

  A \emph{fattened train track}
$\Phi$ on the surface $S$ is a
family of finitely many
  `long' rectangles $e_1$, \dots,
$e_n$  which are foliated
  by arcs parallel to the `short'
sides and which meet only along
arcs (possibly reduced to a
point) contained in their short
sides. In addition, a
fattened train track
$\Phi$ must satisfy the
following:
\begin{enumerate}
  \item[(i)]  each point of the
`short' side of a rectangle also
belongs to  another rectangle,
and each component of the union
of the short sides of all
rectangles is an arc,  as opposed
to a closed curve;

\item [(ii)] each component of the
boundary
$\partial S$ is either disjoint
from
$\Phi$ or contained in it;

  \item[(iii)]  note that the
closure $\overline{S-\Phi}$ of
the complement $S-\Phi $ has a
certain number of `spikes',
corresponding to the points where
at least 3 rectangles meet; we
require that no component of
$\overline{S-\Phi}$ is a disc
with 0, 1 or 2 spikes or an
annulus with no spike.
\end{enumerate} The rectangles
are called the
\emph{edges} of $\Phi$. The
foliations of the edges of $\Phi$
induce a foliation of $\Phi$,
whose leaves are the \emph{ties}
of the fattened train track. The
finitely many ties where several
edges meet are the
\emph{switches} of the
fattened train track
$\Phi$. A tie which is not a
switch is \emph{generic}.
Finally, a point of $\partial
\Phi$ which is contained in three
edges (and consequently is the
tip of a spike of
$\overline{S-\Phi}$) is a
\emph{switch point}.

An $m$--geodesic lamination
$\lambda
$ is \emph{strongly carried} by
the fattened train track $\Phi$ if
it is contained in the interior
of
$\Phi$ in $S$ and if its leaves
are transverse to the ties of
$\Phi$.
  For a fattened train track $\Phi$, let
$\mathcal L^{\mathrm
s}(\Phi)\subset \mathcal L(S)$ be
the subset consisting of those
geodesic laminations which are
strongly carried by $\Phi$.

Given a geodesic lamination
$\lambda$, there are several
constructions which provide a
fattened train track $\Phi$
strongly carrying
$\lambda$; see for instance
\cite{PenHar,
CanEpsGre, Bon10}. In
particular, $\mathcal L(S)$ is
the union of all $\mathcal
L^{\mathrm s}(\Phi)$ as
$\Phi$ ranges over all fattened train
tracks in
$S$. Also, it immediately follows
from definitions that $\mathcal
L^{\mathrm s}(\Phi)$ is open in
$\mathcal L(S)$. Since
$\mathcal L^{\mathrm s}(\Phi)$ is
compact by
Lemma~\ref{thm:L(S)compact}, this
proves:

\begin{lemma}
\label{thm:FinitelyManyTrainTracks}
There exists finitely many fattened train
tracks
$\Phi_1$, $\Phi_2$, \dots,
$\Phi_n$ in
$S$ such that
\begin{equation*}
\mathcal L(S) =
\mathcal L^{\mathrm s}(\Phi_1)
\cup \mathcal L^{\mathrm
s}(\Phi_2)\cup \dots
\mathcal L^{\mathrm s}(\Phi_n) \eqno{\qed}
\end{equation*}
\end{lemma}

Let $a$ be an oriented arc
carried by the fattened train track
$\Phi$, namely such that
$a$ is immersed in the interior of
$\Phi$ and is transverse to its
ties. If $a$ meets the edges
$e_1$, $e_2$,
\dots,
$e_n$ of $\Phi$ in this order, we
will say that $\left\langle e_1,
e_2,
\dots, e_n\right\rangle$ is the
\emph{edge path followed by} the
arc $a$.

If the geodesic lamination
$\lambda$ is strongly carried by
$\Phi$, it
\emph{realizes} the fattened train track
$\left\langle e_1, e_2,
\dots, e_n\right\rangle$ if
$\left\langle e_1, e_2,
\dots, e_n\right\rangle$ is
followed an arc $a$ immersed in a
leaf of $\lambda$.

\section{The Hausdorff dimension
of
$\left(\mathcal L(S), d_{\mathrm
H} \right)$}
\label{sect:MinDim0}

We want to show that the
Hausdorff dimension
$\mathrm{dim_H}
\left( \mathcal L(S), d_{\mathrm
H}\right)$ is equal to 0.  By
Lemma~\ref
{thm:FinitelyManyTrainTracks}, it
suffices to show that, for an
arbitrary fattened train track
$\Phi$,  the subset
$\mathcal  L^{\mathrm s}(\Phi)\subset\mathcal
L(S)$ consisting of those geodesic
laminations which are strongly
carried by
$\Phi$ has Hausdorff dimension 0.

Fix such a fattened train track $\Phi$.
   For a geodesic lamination
$\lambda$ strongly carried by
$\Phi$ and for each  integer
$r\geq 0$, we can consider the set
$\{\gamma_1,
\gamma_2,
\dots\gamma_n\}$ of all edge
paths of length $2r+1$ that are
realized by
$\lambda$. Conversely, given a
finite family of edge paths
$\gamma_1$,
$\gamma_2$, \dots, $\gamma_n$ of
length
$2r+1$, let $P_{\gamma_1 \gamma_2
\dots\gamma_n}\subset \mathcal
L^{\mathrm s}(\Phi)$ consist of those
geodesic laminations such that
the family of edge paths of
length $2r+1$ realized by
$\lambda$ is exactly
$\{\gamma_1, \gamma_2,
\dots\gamma_n\}$, namely such that
$\lambda$ realizes all the
$\gamma_i$ and no other edge path
of length $2r+1$.

\begin{lemma}
\label{thm:DiameterBound} There
exists constants $a$, $b>0$,
depending only on the fattened train track
$\Phi$ and on the negatively
curved metric $m$ of
$S$, such that, for every family
of edge paths
$\gamma_1$,
$\gamma_2$, \dots, $\gamma_n$ of
length
$2r+1$, the diameter of
$P_{\gamma_1
\gamma_2
\dots\gamma_n}\subset \mathcal
L^{\mathrm s}(\Phi)$ for the Hausdorff metric
$d_{\mathrm H}$ is bounded by
$a\mathrm e^{-br}$.
\end{lemma}

\begin{proof} Consider two
geodesic laminations
$\lambda$,
$\lambda'\in P_{\gamma_1
\gamma_2
\dots\gamma_n}$, and their lifts
$\widehat\lambda$,
$\widehat\lambda'$ to the
projective tangent bundle $PT(S)$.

Let $x$ be a point of $\lambda$.
In the leaf of $\lambda$
containing $x$, let $k$ be an
immersed arc  passing through
$x$, and chosen so that each of
the two halves of $k$ delimited
by $x$  crosses exactly
$r+1$ edges (counted with
multiplicities) of
$\Phi$. Let $\gamma$ be the edge
path of length $2r+1$ thus
followed by $k$.

Because $\lambda$ and $\lambda'$
realize the same edge paths of
length
$2r+1$, there is an arc $k'$
immersed in a leaf of $\lambda'$
which follows the same edge path
$\gamma$ as $k$.

Because they follow the same edge
path
$\gamma$, the two $m$--geodesic
arcs
$k$ and $k'$ are homotopic by a
homotopy which moves their end
points by a distance at most
$c_2>0$, where
$c_2$ depends only on the
diameter of the edges of $\Phi$.
The point $x$ separates the arc
$k$ into two arcs whose length is
at least
$c_1 r$, for some constant $c_1>0$
depending on the lengths of the
edges of $\Phi$, suitably defined.
Because the curvature of the
metric $m$ is negative, a Jacobi
field argument then shows that
there exists a point
$x'\in k'$ such that the distance
from
$(x', T_{x'}k')= (x',
T_{x'}\lambda')$ to
$(x, T_{x}k)= (x, T_{x}\lambda)$
in
$PT(S)$ is bounded by
$a\mathrm e^{-br}$, where $a$ and
$b$ depend only on $c_1$, $c_2$
and on a negative upper bound for
the curvature of $m$.

If $\varepsilon = a\mathrm
e^{-br}$, this shows that
$\widehat \lambda$ is contained
in the
$\varepsilon$--neighborhood
$N_\varepsilon\bigl (\widehat
\lambda'\bigr)$ of
$\widehat \lambda'$ in $PT(S)$.

Symmetrically, $\widehat
\lambda'$ is contained in
$N_\varepsilon\bigl (\widehat
\lambda\bigr)$, and it follows
that $d_{\mathrm H}(\lambda,
\lambda')= d_{\mathrm H}(\widehat
\lambda, \widehat \lambda')\leq
\varepsilon = a\mathrm e^{-br}$.
Since this holds for any
$\lambda$,
$\lambda'\in P_{\gamma_1
\gamma_2
\dots\gamma_n}$, this proves that
the diameter of $P_{\gamma_1
\gamma_2
\dots\gamma_n}$ is bounded by
$a\mathrm e^{-br}$.
\end{proof}

Let $\Gamma_r$ denote the set of
edge paths of length $2r+1$ in
$\Phi$. As
$r$ tends to $\infty$, the
cardinal of
$\Gamma_r$ usually grows
exponentially. However, most
elements of $\Gamma_r$ are not
realized by any geodesic
lamination:

\begin{proposition}
\label{thm:PolyBound} The number
of subsets $\{\gamma_1,
\gamma_2, \dots, \gamma_n\}$ of
$\Gamma_r$ for which $P_{\gamma_1
\gamma_2
\dots\gamma_n}$ is non-empty is
bounded by a polynomial function
of $r$.
\end{proposition}
\begin{proof} The argument is in
spirit similar to the core
ingredient of \cite{BirSer} (see
also
\cite[\S10]{Thu86}), which is that
the number of elements of
$\Gamma_r$ that are realized by
geodesic laminations grows at most
polynomially. However, it is
significantly more difficult
because one has no \emph {a
priori} control on the number
$n$ of elements of the subsets
considered.

  Let
$\mathcal P_r$ denote the set of
those subsets
$\{\gamma_1, \gamma_2, \dots,
\gamma_n\}\subset\Gamma_r$ for
which there exists a geodesic
lamination
$\lambda\in P_{\gamma_1, \gamma_2
\dots \gamma_n}$ which
crosses every edge of
$\Phi$. By consideration of the finitely
many fattened train tracks which can be
obtained from
$\Phi$ by removing some edges, it
suffices to show that the number
of elements of $\mathcal P_r$ is
bounded by a polynomial function
of $r$.

Consider an edge path
family $\{\gamma_1,
\gamma_2,
\dots,
\gamma_n\}\in \mathcal P_r$ and a
geodesic lamination
$\lambda \in P_{\gamma_1,
\gamma_2 \dots
\gamma_n}$ crossing every edge of
$\Phi$.
Let $s_1$, $s_2$, \dots, $s_p$ be
the switch points of $\Phi$.

For a given switch point
$s_i$, consider the two leaves
$l$ and
$l'$ of $\lambda$ that are
closest to
$s_i$, namely that can be joined
to
$s_i$ by an arc contained in a
tie and whose interior is
disjoint from
$\lambda$. Note that $l$ and $l'$
exist because $\lambda$ crosses
every edge of
$\Phi$. Then, one of the following
happens:
\begin{enumerate}
\item Starting from the switch
(tie) containing $s_i$, the
leaves $l$ and
$l'$ follow the same edge path of
length $2r+1$.
\item The leaves $l$ and $l'$
follow a common edge path of
length at most
$2r$, and then diverge at some
switch.
\end{enumerate}

In the first case, draw an arc
$z_i$ which is carried by $\Phi$,
is disjoint from $\lambda$,
begins at the switch point $s_i$,
crosses exactly $r$ edges
(counted with multiplicities),
and finally ends on the switch
(tie) at the end of the $r$--th
edge.

In the second case, note that
there is a unique switch point
$s_j$ which separates $l$ and
$l'$ at the switch where they
diverge, and that $l$ and
$l'$ are closest to $s_j$ at that
switch; this again follows from
the fact that
$\lambda$ crosses each edge of
$\Phi$. Then, let
$z_i=z_j$ be an arc which is
carried by
$\Phi$, is disjoint from
$\lambda$ and joins $s_i$ to
$s_j$. Note that
$z_i=z_j$ crosses at most $2r$
edges.

	The arcs $z_i$ are the
$r$--\emph{zippers} associated to
$\lambda$. Those zippers arising
in Case~2 are  called
\emph{switch connections}.

\begin{lemma}
\label{thm:ZippersDefined} The
zippers $z_1$, $z_2$, \dots, $z_p$
are embedded. The zippers $z_i$
and
$z_j$ are disjoint unless they are
equal (and correspond to a switch
connection). Their union
$z^r_\lambda$ is unique up to
isotopy of
$\Phi$ respecting each tie.
\end{lemma}
\begin{proof} We will prove that,
if $Q$ is a component of
$e-\lambda$ where $e$ is an edge
of $\Phi$, the intersection of
$z^r_\lambda$ with $Q$ consists
of at most one arc. The
statements of the lemma
automatically follow from this
property.

For a given $z_i$, consider the
components $Q$ of the above type
that meet $z_i$. By construction,
these form a chain bounded on the
sides by two leaves $l$ and
$l'$ of
$\lambda$ and beginning at the
switch point
$s_i$. It immediately follows that
$z_i$ cannot visit the same $Q$
twice. Similarly, if a $Q$
meeting $z_i$ is also visited by
$z_j$, then the $Q'$ meeting
$z_i\cup z_j$ form a longer
chain, which shows that the two
leaves
$l$ and
$l'$ follow a common edge path of
length at most
$2r$ before diverging at $s_j$; it
follows that
$z_i=z_j$ is a switch connection.
\end{proof}

Let $\mathcal Z_r$ denote the set
of all $r$--zipper families of
this type, considered up to
isotopy of
$\Phi$ respecting its ties.
Namely an element of
$\mathcal Z_r$ is a family $\{
z_1, z_2, \dots, z_p\}$ of
embedded arcs carried by $\Phi$
and of the following two types:
\begin{enumerate}
\item either $z_i$ begins at the
switch point $s_i$, crosses
exactly $r$ edges, and finally
ends on the switch at the end of
the $r$--th edge;
\item or $z_i$ joins the switch
point
$s_i$ to another switch point
$s_j$ and crosses at most $2r$
edges of $\Phi$, in which case
$z_j=z_i$.
\end{enumerate} In addition,
distinct $z_i$ are disjoint. We
identify two such zipper families
when they differ only by an
isotopy of $\Phi$ respecting its
ties.

Lemma~\ref{thm:ZippersDefined}
associates an element
$z^r_\lambda\in \mathcal Z_r$ to
each geodesic lamination
$\lambda$ which is carried by
$\Phi$ and crosses every edge of
$\Phi$.

It is convenient to consider the
set
$\mathcal Q_r$ of all pairs
$\left(\lambda, \{\gamma_1,
\gamma_2,
\dots, \gamma_n \}\right)$ where
$\lambda \in P_{\gamma_1,
\gamma_2,
\dots, \gamma_n}$ crosses every
edge of
$\Phi$.
Lemma~\ref{thm:ZippersDefined}
defines a map $a\co  \mathcal Q_r
\rightarrow
\mathcal Z_r$.

If $2^{\Gamma_r}$ denotes the set
of all subsets of the set
$\Gamma_r$ of all edge paths of
length
$2r+1$, there is also a forgetful
map
$b\co
\mathcal Q_r \rightarrow
2^{\Gamma_r}$ defined by
$b\left(\lambda, \{\gamma_1,
\gamma_2,
\dots, \gamma_n
\}\right)=\{\gamma_1,
\gamma_2,
\dots, \gamma_n \}$. By
construction, its image
$b(\mathcal Q_r)$ is equal to
$\mathcal P_r$.

We now construct a third map
$c\co \mathcal Z_r \rightarrow
2^{\Gamma_r}$.

\begin{lemma}
\label{thm:EdgePathsFromZippers}
Consider an
$r$--zipper family $z\in\mathcal
Z_r$.  Let
$k_1$ and
$k_2$ be two arcs carried by
$\Phi$ which cross $r+1$ edges of
$\Phi$, and which are disjoint
from
$z$. Suppose that
$k_1$ and
$k_2$ start from the same tie $t$
of
$\Phi$ and in the same direction,
and that their starting point is
in the same component of
$t-z$. Then
$k_1$ and $k_2$ follow the same
edge path of length $r+1$.
\end{lemma}
\begin{proof}  Suppose that $k_1$
and $k_2$ diverge at some switch
after following a common edge
path of length $r'<r+1$. Let
$s_i$ be a switch point
separating $k_1$ from
$k_2$ at that switch, and follow
the corresponding component $z_i$
of $z$. If we backtrack from
$s_i$ we see that, for every tie
$t'$ located between $t$ and $s_i$
along $k_1$ and $k_2$, there is a
point of $t'\cap z_i$ which
separates the two points of $k_1$
and $k_2$ located on
$t'$. But for $t'=t$ this would
contradict our hypothesis that the
starting points of $k_1$ and
$k_2$ are in the same component of
$t-z$.
\end{proof}

Lemma~\ref{thm:EdgePathsFromZippers}
defines a map $c\co \mathcal Z_r
\rightarrow 2^{\Gamma_r}$ as
follows. Given $z\in \mathcal
Z_r$, for each edge $e$ of $\Phi$
and each component $f$ of $e-z$,
draw two arcs
$k$ and $k'$ carried by
$\Phi$, starting from the same
point in $f$ but going in opposite
directions, and each crossing
$r+1$ edges before stopping.
Lemma~\ref{thm:EdgePathsFromZippers}
shows that the edge path
$\gamma_f$ followed by $k\cup k'$
depends only on
$f$. Then, $c(z) \in
2^{\Gamma_r}$ is the set of all
edge paths $\gamma_f$ as
$f$ ranges over all components of
$e-z$ and $e$ ranges over all
edges of
$\Phi$.

If $\lambda \in P_{\gamma_1,
\gamma_2, \dots, \gamma_n}$
crosses every edge of $\Phi$, it
follows from the construction of
$z^r_\gamma= a(\lambda,  \{
\gamma_1, \gamma_2, \dots,
\gamma_n\})$ that $\lambda$
crosses every component $f$ of
$e-z_\gamma^r$ for every edge
$e$ of $\Phi$. As a consequence,
$c\circ a=b$, and the  restriction
$c_|\co  a\left(\mathcal
Q_r\right)
\rightarrow b(\mathcal Q_r) =
\mathcal P_r$ is surjective.

The main conclusion of this is
that the cardinal $\#\mathcal
P_r$ of $\mathcal P_r$ is bounded
by
$\#a\left(\mathcal Q_r\right)$,
and therefore by
$\#\mathcal Z_r$ since
$a\left(\mathcal Q_r\right)
\subset \mathcal Z_r$.

\begin{lemma}
\label{thm:Z_rBounded}
If $\Phi$ has $p$
switch points and $q$ edges, then
\begin{equation*}
\# \mathcal P_r \leq
\# \mathcal Z_r \leq 2^p p^{p+q}
r^{p+q}.
\end{equation*} 
\end{lemma}
\begin{proof} We already proved
the first inequality.

Given an $r$--zipper family $z\in
\mathcal Z_r$, we can consider the
number $n_e$ of times it crosses
each edge $e$ of $\Phi$. For each
edge
$e$, the number $n_e$  is bounded
by the total number of edges
crossed by
$z$, counted with multiplicities,
and this number is itself bounded
by $pr$. Consequently, there are
at most $(pr)^q$ possible
assignments for these numbers
$n_e$.

Now, for each edge $e$, $e\cap z$
consists of $n_e$ disjoint
parallel arcs transverse to the
ties. Once we know the $n_e$, we
can read some additional
information from $z$, namely the
location of its end points. More
precisely, for each end point, we
want to  specify the component
$k$ of
$e\cap z$, for some edge $e$ of
$\Phi$, that contains it as well
as which end point of $k$ is the
end point of $z$ considered; note
that there are $\sum_e n_e\leq
pr$ such components of $e\cap z$,
so that we have at most $2pr$
possible choices for each end
point. Some of these end points
are located at the switch points,
and the arc of
$e\cap z$ on which they sit is
completely determined by the
$n_{e'}$. There are at most $p$
of the other end points, which
leaves us with at most  $(2pr)^p$
possible choices.

Once we know that each $e\cap k$
consists of $n_e$ disjoint
parallel arcs transverse to the
ties, as well as the location of
the end points of $z$, then
$z$ is completely determined up to
isotopy of $\Phi$ respecting each
tie. It follows that the number
of elements of $\mathcal Z_r$ is
bounded by $ (pr)^q(2pr)^p= 2^p
p^{p+q}  r^{p+q}$
\end{proof}

Lemma~\ref{thm:Z_rBounded}
concludes the proof of
Proposition~\ref{thm:PolyBound}.
\end{proof}

\begin{remark} The introduction
of the set $\mathcal Q_r$ in the
above proof may at first seem
unnecessarily complicated.
However, it is not hard to find
examples of geodesic laminations
carried by
$\Phi$ which realize exactly the
same edge paths of length $2r+1$
but whose associated $r$--zipper
families are different. In other
words, the map $a\co  \mathcal Q_r
\rightarrow
\mathcal Z_r$ does not
necessarily factor through a
right inverse $\mathcal P_r
\rightarrow
\mathcal Z_r$ for $c$.
\end{remark}

\begin{theorem}
\label{thm:dimd_H=0}

The Hausdorff dimension of $
(\mathcal L(S), d_{\mathrm H})$
is equal to $0$.
\end{theorem}

\begin{proof} By Lemma~\ref
{thm:FinitelyManyTrainTracks}, we
can write $\mathcal L(S)$ as the
union of $\mathcal L^{\mathrm s}(\Phi_i)$ for
finitely many fattened train tracks
$\Phi_i$. It therefore suffices
to show that
$(\mathcal L^{\mathrm s}(\Phi), d_{\mathrm
H})$ has Hausdorff dimension 0
for an arbitrary  fattened train track
$\Phi$.

Let $\Phi$ be such a fattened train track.
By Lemma~\ref{thm:DiameterBound}
and
Proposition~\ref{thm:PolyBound},
there are constants
$a$, $b$, $c$ and
$N>0$ such that, for every integer
$r\geq 1$, we can cover $\mathcal
L^{\mathrm s}(\Phi)$ by at most $cr^N$  sets
of diameter bounded by $a \mathrm
e^{-br}$. Therefore, for every
$\varepsilon>0$, we can cover
$\mathcal L^{\mathrm s}(\Phi)$ by at most
$c\left( 1+\frac1b \log \frac
a\varepsilon
\right)^N$ balls of radius
$\varepsilon$ (by taking $r$ the
smallest integer such that $a
\mathrm e^{-br}< \varepsilon$).
It immediately follows that the
Hausdorff dimension of
$(\mathcal L^{\mathrm s}(\Phi), d_{\mathrm
H})$ is equal to 0.
\end{proof}

For future reference, let us
slightly improve the growth
exponent in
Lemma~\ref{thm:Z_rBounded}.

\begin{lemma}
\label{thm:Z_rBetterBounded} With
the same data as in
Lemma~\ref{thm:Z_rBounded},
\begin{equation*}
\# \mathcal P_r \leq
\# \mathcal Z_r \leq c
r^{9\left\lvert \chi(S)
\right\rvert -1}
\end{equation*} for some constant
$c>0$, where $\chi(S)$ is the
Euler characteristic of $S$.
\end{lemma}

\begin{proof} We first consider
the subset $\mathcal Z_r'$ of
$\mathcal Z_r$ consisting of those
$r$--zipper families which have no
switch connection. We will use the
notation of the proof of
Lemma~\ref{thm:Z_rBounded}.

For $z\in \mathcal Z_r'$, we again
consider the number $n_e$ of
times $z$ crosses the edge $e$ of
$\Phi$. Since
$z$ has no switch connection, the
sum
$\sum_e n_e$ is equal to $pr$. The
numbers $n_e$ are far from being
independent. Indeed, at each
switch
$\sigma$, the sum of the $n_e$
corresponding to the edges $e$
coming in on one side of
$\sigma$ is almost equal to the
sum of the $n_{e'}$ corresponding
to the edges going out on the
other side of
$\sigma$, with a small discrepancy
coming from the end points of $z$
that are located on $\sigma$; in
particular, the discrepancy is
bounded between $-2p$ and $+2p$.

If $z$ and $z'\in\mathcal Z_r'$
have the same discrepancy type at
the switches of $\Phi$, then
their associated families of
crossing numbers $n_e$ and $n_e'$
are such that $n_e' - n_e$ now
satisfy these switch relations
with no discrepancy. The vector
space
$\mathcal W(\Phi)$  of all real
edge weight systems satisfying
the switch relations has
dimension bounded by
$3\left\lvert
\chi(S)\right\rvert$; see for
instance \cite[\S2.1]{PenHar} or
\cite[Theorem~15]{Bon97}. Since
$\sum_e n_e=pr$ for
$z
\in
\mathcal Z_r'$, we conclude that
there exists a constant $c_1$
such that, within a given
discrepancy type at the switches,
the number of possible
assignments of crossing numbers
$n_e$ for $z \in
\mathcal Z_r'$ is bounded by $c_1
r^{3\left\lvert \chi(S)
\right\rvert -1}$ for some
constant $c_1$.

There are at most $(4p+1)^p$
discrepancy types at the switches
of
$\Phi$. The same analysis of end
point locations as in
Lemma~\ref{thm:Z_rBounded} then
gives that $\# \mathcal Z_r' \leq
c_1 (4p+1)^p r^{3\left\lvert
\chi(S)
\right\rvert -1} (2pr)^p$. The
number
$p$ of switch points is bounded by
$6\left\lvert \chi(S)
\right\rvert$ by a counting
argument. It follows that
$\# \mathcal Z_r' \leq c_2
r^{9\left\lvert \chi(S)
\right\rvert -1}$.

Now,  consider the subset
$\mathcal Z_r''$ of
$\mathcal Z_r$ consisting of those
$r$--zipper families which have at
least one switch connection. In
this case, if $z\in \mathcal
Z_r''$ crosses
$n_e$ times the edge $e$, we can
only conclude that $\sum_e n_e
\leq pr$, so that the number of
possible assignments of crossing
numbers $n_e$ for $z\in \mathcal
Z_r''$  is bounded by $c_3
r^{\left\lvert 3\chi(S)
\right\rvert }$ for some constant
$c_3$. However, we now have to
worry about 2 fewer end point
locations, so that
$\# \mathcal Z_r \leq c_3
r^{3\left\lvert \chi(S)
\right\rvert} (2pr)^{p-2}
\leq c_4 r^{9\left\lvert \chi(S)
\right\rvert -2}$.

Since $\mathcal Z_r = \mathcal
Z_r'
\cup \mathcal Z_r''$, this
concludes the proof.
\end{proof}

The estimate of
Lemma~\ref{thm:Z_rBetterBounded}
is
  stronger than that of
Lemma~\ref{thm:Z_rBounded} (where
the growth exponent $p+q$ can be
as large as $15\left\lvert \chi(S)
\right\rvert$), but it is still
quite crude.

\section{The $d_{\log}$ metric}
\label{sect:dlog}

In a metric space $(X,d)$,
consider
\begin{equation*} d_{\log} (x,y)=
\frac1 {\left\lvert \log
\left(
\min\left\{ d(x,y),
\frac14 \right\}
\right)
\right\rvert}
\end{equation*} for $x$, $y\in
X$. In particular,
$d_{\log}(x,y)=1/\left\lvert \log
  d(x,y)
\right\rvert$ when $d(x,y)\leq
\frac14$. As we will see in the
proof of
Proposition~\ref{thm:d_logMetric}
below, the number 4 can be
replaced by any number larger
than 4.

\begin{proposition}
\label{thm:d_logMetric} The above
formula defines a metric
$d_{\log}$ on $X$, which induces
on $X$ the same topology as the
original metric $d$.
\end{proposition}
\begin{proof} The only
non-trivial thing is to prove
that $d_{\log}$ satisfies the
triangle inequality. Because the
function
$f(u)=-1/ \log u$ is increasing
between 0 and 1, it suffices to
show that
$f$ satisfies the inequality
\begin{equation*} f(u+v) \leq
f(u) + f(v)
\end{equation*} for every $u$,
$v\in [0,\frac14]$. For this, we
will compute the minimum of the
function
\begin{equation*} h(u,v)=
f(u)+f(v)-f(u+v)
\end{equation*} over the square
$\left[ 0,
\frac14\right] \times \left[ 0,
\frac14\right] $, and show that
it is equal to 0.

The critical points of $h$ are the
points $(u,v)$ where
$f'(u)=f'(v)=f'(u+v)$. However,
the function
$f'(u)=1/\left(u \log^2 u\right)$
is decreasing on the interval
$\left] 0,
\mathrm e^{-2} \right[$ and
increasing on $\left]
\mathrm e^{-2}, 1\right[$. In
particular, $f'$ is at most
two-to-one on $\left]0,1\right[$.
It follows that a critical point
$(u,v)$ of
$h$ must satisfy $u=v$.  Solving
the equation $f'(u)=f'(2u)$ then
shows that the only critical
point of $h$ in  $\left] 0,
\frac14\right[ \times \left] 0,
\frac14\right[ $ is $\left(
2^{-2-\sqrt 2}, 2^{-2-\sqrt 2}
\right)$, where $h \left(
2^{-2-\sqrt 2}, 2^{-2-\sqrt 2}
\right)=(3-2\sqrt 2)/\log 2 >0$.

On the boundary of the square,
$h(u,0)=h(0,u)=0$ while
$h(u,\frac14)=h(\frac14,u)=g(u)$
with
$g(u)=f(u)+f(\frac14)-f(u+\frac14)$.
 From the graph of $f'$, we see
that $g$ has a unique critical
point, located between $0$ and
$\mathrm e^{-2}$. Graphing $g$,
we conclude that
$g(u)\geq g(0)=g(\frac14)=0$ for
every
$u\in
\left[0,\frac14\right]$.

Therefore, the minimum of $h(u,v)$
over the square $\left[ 0,
\frac14\right] \times \left[ 0,
\frac14\right] $ is equal to 0.
This proves that $f(u+v) \leq
f(u) + f(v)$ for every $u$, $v\in
[0,\frac14]$, and concludes the
proof that
$d_{\log}$ satisfies the triangle
inequality.
\end{proof}

Let us apply this to the case
where
$(X,d)$ is the space $(\mathcal
L(S), d_{\mathrm H})$ of geodesic
laminations. This defines a new
metric $d_{\log}$ on $\mathcal
L(S)$.

This metric $d_{\log}$ is more
intrinsic than $d_{\mathrm H}$.
Recall that two metrics $d$ and
$d'$ on a space $X$ are
\emph{Lipschitz equivalent} is
there exists constants
$C_1$, $C_2>0$ such that
\begin{equation*}
C_1 d(x,y) \leq d'(x,y) \leq C_2
d(x,y)
\end{equation*} for every $x$,
$y\in X$.

\begin{proposition}
\label{thm:d_logLipschitzUnique}
The metric $d_{\log}$ associated
to
$d_{\mathrm H}$ on $\mathcal
L(S)$ is,  up to Lipschitz
equivalence, independent of the
choice of the metric $m$ on
$S$.
\end{proposition}
\begin{proof} This immediately
follows from the fact that
$d_{\mathrm H}$ is well-defined
up to H\"older equivalence
(Lemma~\ref{lem:HolderInv}).
\end{proof}

In particular, the Hausdorff
dimension of
$(\mathcal L(S), d_{\log})$ is
now well-defined, independently
of the choice of a metric on $S$.

We have the following
estimate for this Hausdorff
dimension
$\mathrm{dim_H}(\mathcal L(S),
d_{\log})$. This estimate is far
from sharp, but at least proves
that the  dimension is
  strictly between $0$ and
$\infty$.

\begin{proposition}
\label{thm:dimRoughBounds} If the
surface
$S$ is neither the three times
punctured sphere, nor the twice
punctured projective plane, nor
the once punctured Klein bottle,
\begin{equation*} 2 \leq
\mathrm{dim_H} (\mathcal L(S),
d_{\log})
  \leq 9\left\lvert\chi(S)
\right\rvert-1
\end{equation*} where $\chi(S)<0$
is the Euler characteristic of
$S$.
\end{proposition}

As indicated in the introduction,
$\mathcal L(S)$ is finite when
$S$ is  the three times punctured
sphere,  the twice punctured
projective plane, or the once
punctured Klein bottle; see for
instance
\cite[Proposition~16]{BonZhu3}
for a proof. In particular,
$\mathcal L(S)$ has Hausdorff
dimension equal to 0 for these
exceptional surfaces.

\begin{proof}[Proof of
Proposition~\ref{thm:dimRoughBounds}]
The second inequality is a
by-product of the proof of
Theorem~\ref{thm:dimd_H=0}.
Indeed, for every integer $r\geq
1$,
  Lemmas~\ref{thm:DiameterBound}
and
\ref{thm:Z_rBetterBounded} show
that  we can cover
$\mathcal L(S)$ by at most $c
r^{9\left\lvert\chi(S)
\right\rvert-1}$ sets whose
$d_{\mathrm H}$--diameter is
bounded by $a\mathrm e^{-br}$,
for some constants $a$, $b$,
$c>0$. By definition of
$d_{\log}$, the
$d_{\log}$--diameter of these
sets is bounded by $a'r^{-1}$ for
some constant
$a'$ depending only on $a$ and
$b$. It immediately follows that
$\mathrm{dim_H}
  (\mathcal L(S), d_{\log}) \leq
9\left\lvert\chi(S)
\right\rvert-1$.

The first inequality is a
consequence of \cite{BonZhu3},
which itself uses \S\ref
{sect:TrainTracks} of the current
paper. Indeed, if
$S$ is neither the three times
punctured sphere, nor the twice
punctured projective plane, nor
the once punctured Klein bottle,
then $S$ contains a simple closed
geodesic
$\gamma$ such that one component
$T$ of the surface obtained by
splitting
$S$ open along
$\gamma$ is either a once
punctured torus or a  four times
punctured sphere. Then, the set
$\mathcal L_0(T)$, consisting  of
those geodesic laminations which
are contained in the interior of
$T$, is a natural subspace of
$\mathcal L(S)$. Since we
show in \cite{BonZhu3} that
$\mathrm {dim_H} (\mathcal L_0(T),
d_{\log}) =2$, this proves the
first inequality.
\end{proof}

\section{The combinatorial
distance $d_{\Theta}$}
\label{sect:TrainTracks}

Having used fattened tracks so
far, it is now convenient to
switch to (unfattened) train
tracks.
A \emph{train
track} on  the surface $S$
  is a graph $\Theta $ contained in
the interior of $S$ which
consists of finitely many
vertices, also called
\emph{switches}, and of finitely
many  edges joining them such
that:
\begin{enumerate}
\item    the edges of $\Theta $
are  differentiable arcs whose
interiors are embedded  and
pairwise disjoint (the two end
points of an  edge may coincide);
\item    at each switch $s$ of
$\Theta
$,  the edges of $\Theta $ that contain $s$
are all  tangent to the same line $L_{s}$ in
the tangent  space $T_{s}S$ and, for each of
the two  directions of $L_{s}$, there is at
least one edge  which is tangent to that
direction;
\item    observe that the
complement
$S-\Theta $ has a certain number of spikes,
each  leading to a switch $s$ and locally
delimited by  two edges that are tangent to
the same direction  at $s$; we require that no
component of $S-\Theta
$ is a disc with 0, 1 or 2 spikes or an open
annulus with no spike .
\end{enumerate}

A curve $c$
\emph{carried} by the train track
$\Theta$ is a differentiable
immersed curve in $c\co I
\rightarrow S$ whose image is
contained in
$\Theta$, where $I$ is an
interval in $\mathbb R$. Such a
curve is
\emph{bi-infinite} if its
restriction to each component of
$I-\{x\}$ has infinite length,
where $x$ is an arbitrary
point in the interior of $I$.

In the definition of train
tracks, the third condition is
particularly crucial, as it has
the following global corollaries
(see for instance
\cite[\S1.8]{Bon10} for proofs).

\begin{lemma}
\label{lem:CurvesTrainTracks}
Let $\Theta$ be a train track in
the surface $S$, and let
$\widetilde \Theta$ be its
preimage in the universal
covering $\widetilde S$ of $S$.
Then,
\begin{enumerate}
\item any curve carried by
$\widetilde \Theta$ is embedded;
\item if two curves carried by
$\widetilde \Theta$
coincide for a while and then
diverge at some switch, they never
meet again;
\item any bi-infinite curve which
is carried by $\widetilde\Theta$
is quasi-geodesic in $\widetilde
S$.\qed
\end{enumerate}
\end{lemma}

Consider a bi-infinite curve $c$
carried by the train track
$\Theta$. By
Lemma~\ref{lem:CurvesTrainTracks}
  and by the classical
property of quasi-geodesics in
negatively curved manifolds,
there exists  a unique
geodesic $g$ which, after a
possible reparametrization of
$c$, is
homotopic to $c$ by a homotopy
which moves points by a uniformly
bounded distance. In this
situation, we will say that the
geodesic $g$ is
\emph{tracked by the curve} $c$
carried by $\Theta$.

The geodesic lamination $\lambda$
is
\emph{weakly carried} by the train
track $\Theta$ if every leaf of
$\lambda$ is tracked by a
bi-infinite curve carried by
$\Theta$. Let $\mathcal L^{\mathrm
w} (\Theta)$ denote the subset of
$\mathcal L(S)$ consisting of
those geodesic laminations which
are weakly carried by $\Theta$.

As the name indicates, a
fattened train track
$\Phi$ is a certain thickening of
a train track $\Theta$. Every
geodesic lamination which is
strongly carried by the
fattened train track $\Phi$ is
weakly carried by the train track
$\Theta$, so that $\mathcal
L^{\mathrm s}(\Phi) \subset
\mathcal L^{\mathrm w}(\Theta)$.
Each point of view has its own
technical advantages, and the
framework of train tracks is
somewhat better adapted to the
context of the remainder of this
paper. In particular,
unlike $\mathcal L^{\mathrm
s}(\Phi)$,
$\mathcal L^{\mathrm w}(\Theta)$
does not depend on the choice of
a negatively curved metric on
$S$.

By Lemma~\ref
{thm:FinitelyManyTrainTracks},
there exists  finitely many train
tracks
$\Theta_1$,
\dots,
$\Theta_n$ such that
$\mathcal L(S)$ can be written as
the union of the
$\mathcal L^{\mathrm w}
(\Theta_i)$. See
\cite{PenHar} for an explicit such
family of $\Phi_i$. Our goal is to
understand the metrics
$d_{\mathrm H}$ and
$d_{\log}$ on each of these
$\mathcal L^{\mathrm w}
(\Theta_i)$.

A bi-infinite curve carried by
the train track $\Theta$ crosses
a succession of oriented edges
\dots, $e_{-n}$, \dots, $e_0$,
$e_1$, \dots, $e_n$, \dots, in
this order. The list
$\left\langle
\dots, e_{-n}, \dots, e_0,
e_1, \dots, e_n, \dots
\right\rangle$ is the
\emph{bi-infinite edge path}
followed by the curve $c$. If the
geodesic
$g$ is tracked by a curve $c$ on
the train track
$\Theta$, the bi-infinite train
track followed by $c$ in $\Theta$
depends uniquely on
$g$; see \cite[\S1.8]{Bon10}. A
finite length edge path $\gamma$
is
\emph{realized} by the geodesic
lamination $\lambda$ if it is
contained in the bi-infinite edge
path associated to a leaf of
$\gamma$.

If $\lambda$ is weakly carried by
$\Theta$ and if $r\geq 1$ is an
integer, we can consider the set
of all the edge paths of length
$r$ of $\Theta$ which are realized
by
$\lambda$. Given two geodesic
laminations $\lambda$,
$\lambda'\in
\mathcal L^{\mathrm w} (\Theta)$, we
then define
\begin{equation*}
\textstyle
d_\Theta(\lambda,
\lambda')=
\min \left\{ \frac1{r+1};
\lambda \text{ and } \lambda'
\text{ realize the same edge
paths of length } r
\right\}.
\end{equation*}

This $d_\Theta(\lambda, \lambda')$
is the
\emph{combinatorial distance}
between $\lambda$ and $\lambda'$
in  $\mathcal L^{\mathrm
w}(\Theta)$.

\begin{lemma}
\label{thm:d_PhiUltrametric} The
function $d_\Theta$ defines an
ultrametric on the set $\mathcal
L^{\mathrm w} (\Theta)$.
\end{lemma}

Recall that an \emph{ultrametric}
on a set $X$ is a metric $d$
where the triangle inequality is
replaced by the stronger
condition that $d(x,z)
\leq \max \{ d(x,y), d(y,z)\}$ for
every $x$, $y$, $z\in X$.

\begin{proof}[Proof of
Lemma~\ref{thm:d_PhiUltrametric}]
The only non-trivial thing to
prove is that $d_\Theta(\lambda,
\lambda')$ can be $0$ only if
$\lambda=\lambda'$. (In
particular, note that the
ultrametric inequality is
completely tautological.)

Suppose that $d_\Theta(\lambda,
\lambda')=0$, namely that
$\lambda$ and $\lambda'$ realize
exactly the same edge paths.

Let $g$ be a leaf of $\lambda$.
If $\gamma$ is an edge path
contained in the bi-infinite edge
path realized by $g$, there
consequently exists a leaf
$g_\gamma'$ of
$\lambda'$ whose associated
bi-infinite edge path contains
$\gamma$. Taking increasing larger
edge paths $\gamma$ and passing
to a converging subsequence of
leaves of
$\lambda$, we conclude that there
exists a leaf $g'$ of $\lambda'$
which follows exactly the same
bi-infinite edge path as $g$. In
this case, the two geodesics $g$
and $g'$ are homotopic by a
homotopy which moves points by a
uniformly bounded amount. Since
the curvature of the metric $m$
is negative, this implies that
$g'=g$.

Since this property holds for
every leaf $g$ of $\lambda$, it
follows that
$\lambda\subset\lambda'$, and
therefore that $\lambda =
\lambda'$ by symmetry.
\end{proof}

\section{Hausdorff
distance and
combinatorial distance}
This section is devoted to
proving Theorem~\ref
{thm:HausCombEstimate1} of the
Introduction, which we restate
here as:
\begin{theorem}
\label{thm:HausCombEstimate2}
  Let $\Theta$ be a
train track on the
negatively curved surface
$S$. Then, there exists two
constants $c_1$, $c_2>0$
such that, for every geodesic
laminations
$\lambda$, $\lambda'$ which are
weakly carried by $\Theta$,
\begin{equation*}
\mathrm e^{-c_1/d_{\Theta}
\left( \lambda,
\lambda' \right)}
\leq
d_{\mathrm H} \left( \lambda,
\lambda' \right)
\leq  \mathrm e^{-c_2 /d_{\Theta}
\left( \lambda,
\lambda' \right)}
\end{equation*}
\end{theorem}

We will split the proof of
Theorem~\ref
{thm:HausCombEstimate1} into two
lemmas, each devoted to one of
the two inequalities. We begin
with the easier one, namely the
second inequality.
\begin{lemma}
\label{thm:d_log<d_Theta}
There
exists  constants $c$, $c'>0$
such that
\begin{equation*} d_{\mathrm H}
(\lambda, \lambda')
\leq c \thinspace \mathrm e^{-c'/
d_\Theta (\lambda,
\lambda')}
\end{equation*} for every
$\lambda$, $\lambda' \in
\mathcal L^{\mathrm w} (\Theta)$.
\end{lemma}

\begin{proof}
The proof
is in spirit very similar to that
of Lemma~\ref{thm:DiameterBound},
except that we now have to worry
about quasi-geodesics.
  We can assume
$\lambda\not= \lambda'$, since
otherwise there is nothing to
prove, and  we can restrict
attention to the case where
$d_\Theta (\lambda,
\lambda')\leq\frac12$ since
$d_{\mathrm H}$ is bounded.
Consider the integer
$r\geq0$ such that
$d_\Theta (\lambda,
\lambda')=\frac1{2r+2}$ or
$\frac1{2r+3}$. In particular,
$\lambda$ and
$\lambda'$  realize the same edge
paths of length $2r+1$.

Consider a point
$x\in\lambda$. Since  $\lambda \in
\mathcal L^{\mathrm w} (\Theta)$ is
weakly carried by $\Theta$,  the
leaf
$g$ of
$\lambda$ passing through $x$ is
tracked by a bi-infinite curve $c$
carried by
$\Theta$. In particular, there is a
homotopy from $g$ to
$c$ which moves points by a
distance bounded by a constant
$c_1>0$ depending only on $\Theta$
and on the metric $m$.  Let $y\in
c$ be the image of $x\in g$ under
this homotopy. Let $a$ be an arc
contained in $c$, containing $y$
and crossing exactly $r+1$ edges
in each direction when starting
from
$y$, and let $\gamma$ be the edge
path of length $2r+1$ that is
followed by $a$.  Let $b\subset
g$ be the image of $a$ by the
homotopy from $c$ to $g$.

Because $d_\Theta (\lambda,
\lambda')\leq\frac 1{2r+2}$,
$\lambda'$ also realizes the edge
path
$\gamma$. Namely, there is a leaf
$g'$ of $\lambda'$ which is
tracked by a bi-infinite curve
$c'$  carried by $\Theta$, and an
arc
$a'\subset c'$ which follows
$\gamma$.

The two arcs $a$ and $a'$ are
homotopic by a homotopy which
moves points by a distance $\leq
c_2$, where $c_2$ depends only on
the diameters of the edges of
$\Theta$, and $a'$ is also
homotopic to an arc $b'\subset
g'$ by a homotopy which moves
points by a distance $\leq c_1$.

Now, the two geodesic arcs
$b\subset g$ and $b'\subset g'$
are homotopic by a homotopy which
moves their end points by at most
$ 2c_1+c_2$. In addition, by
quasi-geodesicity of $c$, the
point
$x$ is at distance at least $c_3
r$ from the two end points of $b$,
for some constant $c_3>0$. Since
the curvature of
$m$ is negative, it follows that
there exists an $x'\in b'$  such
that
$(x, T_xb)$ is at distance at most
$c_4 \mathrm{e}^{-c_5r}$ from
$(x', T_{x'}b')$ in the
projective tangent bungle
$PT(S)$, for some constants
$c_4$, $c_5>0$ depending only on
$\Theta$ and on the metric $m$.

Since this property holds for
every
$x\in \lambda$, we conclude that
the lift $\widehat\lambda$ of
$\lambda$ in
$PT(S)$ is contained in the $c_4
\mathrm{e}^{-c_5r}$--neighborhood
of the lift $\widehat\lambda'$ of
$\lambda'$. Symmetrically,
$\widehat\lambda$ of $\lambda$ in
$PT(S)$ is also contained in the
$c_4
\mathrm{e}^{-c_5r}$--neighborhood
of
  $\widehat\lambda'$, and it
follows that
\begin{equation*} d_{\mathrm H}
(\lambda, \lambda')= d_{\mathrm
H} (\widehat\lambda,
\widehat\lambda') \leq c_4
\mathrm{e}^{-c_5r}
\leq c_6
\mathrm{e}^{-c_5/2
d_\Theta(\lambda, \lambda')}
\end{equation*}
for $c_6= c_4 \mathrm e^{3c_5/2}$
This proves the
inequality of
Lemma~\ref{thm:d_log<d_Theta}.
\end{proof}

We now prove the first equality
of Theorem~\ref
{thm:HausCombEstimate1}. Its
proof is somewhat more elaborate
than that of the other
inequality, and has a more
topological flavor.

\begin{lemma}
\label{thm:d_Teta<d_log} There
exists  constants $c$, $c'>0$
such that
\begin{equation*} d_{\mathrm
H} (\lambda, \lambda')
\geq c\thinspace \mathrm e^{-c'/
d_{\Theta} (\lambda,
\lambda')}
\end{equation*} for every
$\lambda$, $\lambda' \in
\mathcal L^{\mathrm w} (\Theta)$.
\end{lemma}

\begin{proof}
  It is convenient to arrange, by
adding a few edges to
$\Theta$ if necessary, that
$\Theta$ is maximal for
inclusion in the set of all train
tracks. This is equivalent to the
property that the complement
$S-\Theta$ consists of finitely
many triangles and finitely
semi-open annuli (each
containing a component of
$\partial S$).

  To simplify the
notation, let the integer
$r\geq0$ be such that
$d_\Theta (\lambda, \lambda')
=\frac1{2r}$ or $\frac1{2r+1}$
(assuming $\lambda\not=\lambda'$
without loss of generality), and
let
$\delta=2d_{\mathrm H} (\lambda,
\lambda')$.

Since $d_\Theta (\lambda, \lambda')
>\frac1{2r+2}$ and exchanging the
r\^oles of $\lambda$ and
$\lambda'$ if necessary, there
exists an edge path
$\gamma$ of length $2r+1$ which is
realized by $\lambda$ and not by
$\lambda'$.

By definition, the fact that
$\lambda\in \mathcal L^{\mathrm
w}(\Theta)$ realizes $\gamma$ means
that there exists a leaf $g$
which is tracked by a bi-infinite
curve $c$ carried by
$\Theta$, and an arc $a\subset c$
which  follows the edge path
$\gamma$. Pick a point $y\in c$
in the part of $a$ corresponding
to the $(r+1)$--th edge which it
crosses, namely to the central
edge of $\gamma$. By
quasi-geodesicity, there is a
homotopy from $c$ to $g$ which
moves points by a distance
bounded by
$c_1$. Let $x\in g$ be the image
of
$y$ under this homotopy.

Since $d_{\mathrm H} (\lambda,
\lambda')<\delta$, there is an
$x'\in\lambda'$, located in a leaf
$g'$ of $\lambda'$, such that the
distance from $(x, T_xg)$ to
$(x',T_{x'}g')$ in $PT(S)$ is
bounded by $\delta$. The leaf
$g'$ is tracked by $c'$ carried by
$\Theta$, and there is a homotopy
from $c'$ to $g'$ which moves
points by a distance bounded by
$c_1$.

At this point, it is convenient to
lift the situation to the
universal covering $\widetilde S$
of $S$. Let
$\widetilde\Theta$, $\widetilde
\lambda$ and $\widetilde
\lambda'$ be the respective
preimages of $\Theta$,
$\lambda$ and $\lambda'$. Lift
$x$ to a point $\widetilde x$,
contained in a leaf $\widetilde
g$ of
$\widetilde\lambda$. Lift $x'$ to
a point
$\widetilde x'$, contained in the
leaf
$\widetilde g'$ of $\widetilde
\lambda'$, such that the distance
from
$(\widetilde x, T_{\widetilde x}
\widetilde g)$ to $(\widetilde x',
T_{\widetilde x'}
\widetilde g')$ in
$PT\bigl(\widetilde S\bigr)$ is
bounded by $\delta$. Lift $c$ and
$c'$ to bi-infinite curves
$\widetilde c$ and $\widetilde
c'$ carried by
$\widetilde \Theta$, in such a way
that they are respectively
homotopic to
$\widetilde g$ and $\widetilde
g'$ by a homotopy which moves
points by a distance bounded by
$c_1$. Finally, the  homotopy
from $\widetilde g$ to
$\widetilde c$ specifies preferred
lifts
$\widetilde  y$,
$\widetilde  a$ of $y$, and $a$,
as well as a lift of $\gamma$ to
an edge path $\widetilde \gamma$
of
$\widetilde \Theta$.

Note that $\widetilde g$ and
$\widetilde c$ have the same end
points on the boundary at infinity
$\partial_\infty \widetilde S$.
Similarly, $\widetilde g'$ and
$\widetilde c'$ have the same end
points in
$\partial_\infty \widetilde S$.
Since the distance from
$(\widetilde x, T_{\widetilde x}
\widetilde g)$ to $(\widetilde x',
T_{\widetilde x'}
\widetilde g')$ in
$PT\bigl(\widetilde S\bigr)$ is
bounded by $\delta$, we conclude
that each end point of
$\widetilde  c$ is seen from
$\widetilde x$ within an angle of
$\leq c_7\delta$ from an end
point of $\widetilde c'$, where
the constant $c_7$ depends only
on $c_1$ and on the curvature of
$m$.

We now split the argument into
three cases.

\begin{figure}[ht!]
\begin{center}
\includegraphics{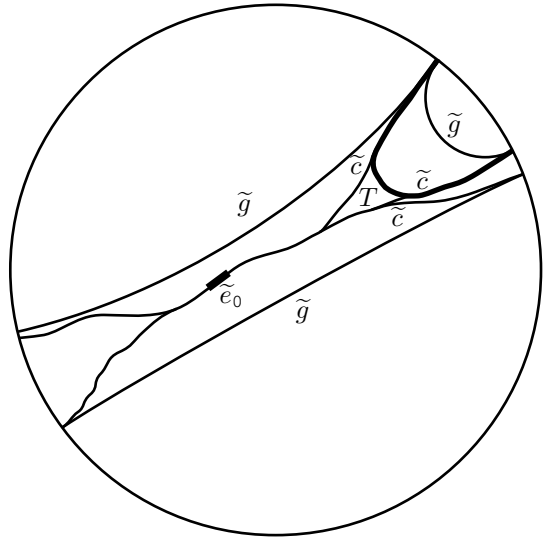}
\end{center}
\caption{Diverging curves carried
by the train track $\widetilde
\Theta$}
\label{pict:DivergCurves}
\end{figure}

\medskip\noindent\textbf{Case 1}\qua
\textsl{The curve $\widetilde c'$
crosses the edge $\widetilde e_0$
of
$\widetilde \Theta$ that contains
the point $y\in\widetilde c$. }

By choice of $\gamma$, the curve
$\widetilde c'$ does not realize
the edge path $\widetilde \gamma$.
Consequently,  we can follow
$\widetilde c$ and
$\widetilde c'$ in one direction
from
$\widetilde e_0$ until they
diverge at some switch
$\widetilde s$ of
$\widetilde
\Theta$ after crossing at most $r+1$
edges of $\widetilde \Theta$.
Orient
$\widetilde c$ and $\widetilde
c'$ in the direction from
$\widetilde e_0$ to $\widetilde
s$.

Since $\widetilde c$ and
$\widetilde c'$ diverge at the
switch $\widetilde s$, there is a
spike of
$\widetilde S-\widetilde
\Theta$ which separates
$\widetilde c$
from $\widetilde c'$  at
$\widetilde s$.  Let
$T$ be the component of
$\widetilde S-\widetilde \Theta$
that contains that spike. By the
maximality hypothesis on $\Theta$,
either $T$ is a triangle or it is
the preimage of a semi-open
annulus component of $S-\Theta$.

We first analyze the case where
$T$ is a triangle.  Consider the
side of
$T$ that does not contain
$\widetilde s$. Extend this side
of $T$
to a bi-infinite curve
$\widetilde c''$ which is carried
by
$\widetilde \Theta$. By the
second conclusion of Lemma~\ref
{lem:CurvesTrainTracks},
either $\widetilde c''$ remains
disjoint from $\widetilde c$, or
it meets $\widetilde c$ at one
spike of $T$; in the second case,
we can therefore assume that
$\widetilde c''$ coincides with
$\widetilde c$ after leaving the
closure of
$T$. Similarly, we can assume
that either $\widetilde c''$ is
disjoint from $\widetilde c'$ or
it coincides with $\widetilde c'$
after leaving the closure $T$. See
Figure~\ref{pict:DivergCurves}.
Again,
$\widetilde c''$ is
quasi-geodesic, and there is a
geodesic
$\widetilde g''$ which is
homotopic to
$\widetilde c''$ by a homotopy
moving points by a distance
bounded by $c_1$.

By construction, the distance from
$\widetilde y$ to the geodesic
$\widetilde g''$ is bounded by
$c_8 r$ for a constant $c_8$
depending only on $c_1$, the
length of the edges of
$\Theta$, and the diameter of the
components of $S-\Theta$.  It
follows that the angle under
which its end points are seen
from $\widetilde y$ is at least
$\mathrm e^{-c_9 r}$, where
$c_9$ depends on $c_8$ and on a
negative  bound for the curvature
of
$m$. By construction, these end
points separate the positive end
points of
$\widetilde c$ and $\widetilde
c'$. We conclude that these
positive end points are seen
from $\widetilde x$ under an
angle larger than $e^{-c_9 r}$.

On the other hand, splicing one
half of $\widetilde c$ to one
half of
$\widetilde c'$ in $\widetilde
e_0$, we can construct a curve
$\widetilde c'''$ which is
carried by
$\widetilde\Theta$, which crosses
$\widetilde e_0$, and which goes
from the negative end point of
$\widetilde c'$ to the positive
end point of
$\widetilde c$. By
quasi-geodesicity of $\widetilde
c'''$, it follows that the angle
under which these two end points
are seen from $\widetilde x$ is
bounded from below by a constant
$c_{10}>0$.

	However, we had concluded earlier
  that each end point of
$\widetilde  c$ is seen from
$\widetilde x$ within an angle of
$\leq c_7\delta$ from at least one
end point of $\widetilde c'$. It
follows that $ e^{-c_9 r} \leq c_7
\delta$ for $r$ sufficiently
large.

There remains to consider the
case where $T$ is the preimage of a semi-open
annulus component of $S-\Theta$.
Namely, $T$ is an infinite strip
containing a component
$\widetilde g''$ of
$\partial \widetilde S$.  The
geodesic
$\widetilde g''$ separates the
two positive end points of
$\widetilde c$ and $\widetilde
c'$, and its distance from
$\widetilde y$ is again  bounded
by
$c_8 r$ for a constant $c_8$
depending only on the
length of the edges of
$\Theta$ and on the diameter of
the components of $S-\Theta$. The
same argument as before then
shows that $ e^{-c_9 r} \leq c_7
\delta$ for $r$ sufficiently
large in this case as well.

\medskip\noindent\textbf{Case 2}\qua
\textsl{The curve
$\widetilde c'$ does not cross
the edge $\widetilde e_0$
containing $y\in\widetilde c$,
but meets $\widetilde c$ at some
edge $\widetilde e_1 \not =
\widetilde e_0$.}

Orient $\widetilde c$ in the
direction from $\widetilde e_1$
to $\widetilde e_0$, and orient
$\widetilde c'$ so that its
orientation coincides with that
of $\widetilde c$ on $\widetilde
e_1$.

The curve $\widetilde c'$
diverges from $\widetilde c$ at
some point between $\widetilde
e_1$ and
$\widetilde e_0$. In particular,
the direction in which it diverges
specifies a side of $\widetilde
c$. Among the two components of
$\widetilde S-\widetilde \Theta$
which are adjacent to $\widetilde
e_0$, let
$T$ be the one which is on this
preferred side of $\widetilde c$.
Again, either $T$ is a triangle or
  it is
the preimage of a semi-open
annulus component of $S-\Theta$.

First consider the case where $T$
is a triangle.
Among the two sides of
$T$ which do not touch $\widetilde
e_0$, the curve $\widetilde c'$
may meet one of them in
$\widetilde
\Theta$, but not both. Pick a side
of
$T$ which is not followed by
$\widetilde c'$ and
extend it to a bi-infinite curve
$\widetilde c''$ which is carried
by
$\widetilde \Theta$. As in
Case~1, we can arrange that
$\widetilde c''$ either is
disjoint from
$\widetilde c$ (resp. $\widetilde
c'$), or coincides with this
curve after it meets it.

This time the geodesic
$\widetilde g''$ homotopic to
$\widetilde c''$ passes at
uniformly bounded distance from
$\widetilde x$. Since it
separates the positive end points
of
$\widetilde c$ and
$\widetilde c'$, we conclude
that the angle under which these
two end points are seen from
$\widetilde x$ is uniformly
bounded from below by $c_{10}>0$.

As in Case~1, the angle under
which the positive end point of
$\widetilde c$ and the negative
end point
$\widetilde c'$s are seen from
$\widetilde x$ is bounded from
below by $c_{10}>0$.

We conclude in this case that
$c_{10}<c_7 \delta$, and therefore
that $ e^{-c_9 r} \leq c_7
\delta$ for $r$ sufficiently
large.

The argument is similar when $T$
is the preimage of a semi-open
annulus component of $S-\Theta$,
by considering the geodesic
component $\widetilde g''$ of
$\partial S$ contained in $T$.

\medskip\noindent\textbf{Case 3}\qua
\textsl{ The curves $\widetilde c$
and $\widetilde c'$ never cross
the same edge.}

The same argument as in Case~2,
using the component $T$ of
$\widetilde S-\widetilde \Theta$
adjacent to
$\widetilde e_0$ and on the same
side as $\widetilde c'$, gives
that
$c_{10}<c_7 \delta$, and therefore
that $ e^{-c_9 r} \leq c_7
\delta$ for $r$ sufficiently
large.

\medskip

Therefore, in all cases,
$d_{\mathrm H} (\lambda,
\lambda') = \frac12 \delta \geq
\frac12 c_{7}^{-1} e^{-c_{9} r}
\geq
\frac12 c_{7}^{-1} e^{-c_{9}
/2d_\Theta (\lambda, \lambda')}$
for $r$ sufficiently large. This
concludes the proof of
Lemma~\ref{thm:d_Teta<d_log}.
\end{proof}

The combination of Lemmas~\ref
{thm:d_log<d_Theta}
and \ref{thm:d_Teta<d_log} proves
Theorem~\ref
{thm:HausCombEstimate2}.

\end{document}